\documentclass[11pt, reqno]{amsart}
\usepackage{amssymb,latexsym,amsmath,amsfonts}

\numberwithin{equation}{section}

\theoremstyle{definition}
\newtheorem{definition}{Definition}[section]
\newtheorem{example}[definition]{Example}
\newtheorem{custom}[definition]{}
\theoremstyle{remark}
\newtheorem{remark}[definition]{Remark}

\theoremstyle{plain}
\newtheorem{theorem}[definition]{Theorem}
\newtheorem{corollary}[definition]{Corollary}
\newtheorem{lemma}[definition]{Lemma}
\newtheorem{proposition}[definition]{Proposition}

\newtheorem{result}[definition]{Result}

\newcommand{\eps}{\varepsilon}

\newcommand{\zt}{\zeta}

\def \disk{\mathbb{D}}
\def \hol{\mathcal{O}} 
\def \C#1{\mathbb{C}^{#1}}

\def \cplx{\mathbb{\mathbb{C}}}

\newcommand{\bdy}{\partial}

\newcommand{\bcdot}{\boldsymbol{\cdot}}

\begin{document}

\title{Proper holomorphic mappings 
between hyperbolic product manifolds}

\author{Jaikrishnan Janardhanan}
\address{Department of Mathematics, Indian Institute of Science, Bangalore -- 560012}
\email{jaikrishnan@math.iisc.ernet.in}

\thanks{This work is supported by a UGC Centre for Advanced Study grant
and by a scholarship from the IISc.}


\begin{abstract}
  We prove a result on the structure of finite proper 
  holomorphic mappings between complex manifolds that are products of hyperbolic 
  Riemann surfaces. While an important special case of our result follows from
  the ideas developed by Remmert and Stein, the proof of the full result relies
  on the interplay of the latter ideas and a finiteness theorem for Riemann
  surfaces.
\end{abstract} 
\maketitle

\section{Introduction}
A natural problem in the study of holomorphic mappings is 
to give a description of the automorphism group or the set of all proper 
holomorphic self-maps of a given complex space. Of course, this problem is
intractable in general even when posed for arbitrary domains in $\C{n}$. One of the
simplest classes of domains in $\C{n}$ is the class of products of planar domains. 
The following theorem due to Remmert and Stein gives the precise structure of a
proper holomorphic map between certain product domains. 

\begin{result}[Remmert--Stein \cite{remmert1960eha}]\label{RES:RS}
  Let $D = D_1\times \cdots \times D_n$ and $W = W_1\times \cdots \times W_n$
  be products of planar domains such that, for each $j = 1, \ldots, n$,
  $\cplx \setminus D_j$ has non-empty interior and $W_j \subset \cplx$ is
  bounded. Let $f=(f_1,\ldots,f_n)$  be a proper holomorphic map from 
  $D$ to $W$. Then, each $f_j,\ j = 1, \ldots, n$, is of the form 
  $f_j(z_{p(j)})$, where $p$ is a permutation of $\{1,\ldots,n\}$. 
\end{result}

\begin{remark}\label{RMK:RADO}
  The proof of this result in the case $n=2$ was given by Remmert and Stein
  (Satz $12$ in \cite{remmert1960eha}). Their proof uses Rado's theorem. 
  The proof of the general case requires a generalization of Rado's theorem,
  but all other aspects of Remmert and Stein's proof
  remain unchanged: see \cite[pp.~71--78]{narasimhan1971scv}, for instance.
\end{remark}

The methods of Remmert--Stein cannot, as Remark~\ref{RMK:RADO} suggests, be
applied when the planar factors in Result~\ref{RES:RS} are replaced by compact
Riemann surfaces. Yet, the phenomenon exhibited by Result~\ref{RES:RS} occurs
in some other settings as well.  Another result on mappings between product
spaces is the following theorem of Peters which generalizes a well-known result by
Cartan \cite{cartan1936domainesBornes}:

\begin{result}[Peters \cite{peters1974Produkte}]\label{RES:AUT}
  Let $X$ and $Y$ be hyperbolic complex spaces. Then the natural injection
  $\text{Aut}(X) \times \text{Aut}(Y) \to \text{Aut}(X \times Y)$ induces an
  isomorphism
  \begin{equation*}
	\text{Aut}^0(X) \times \text{Aut}^0(Y) \cong \text{Aut}^0(X \times Y).
  \end{equation*}
\end{result}
\noindent Here $\text{Aut}^0(X)$ denotes the connected component of the identity
element of $\text{Aut}(X)$.

\smallskip
To the best of our knowledge, there is no analogue of the above result for proper 
holomorphic maps in the literature, except for Result \ref{RES:RS} (and a small
technical improvement thereof in \cite[p.~77]{narasimhan1971scv}). It is natural
to ask whether a result similar to Result~\ref{RES:RS}, with the planar domains therein 
replaced by hyperbolic Riemann surfaces, holds true. However, it is far from clear whether the 
methods seen in the proofs of {\em either} of the above results are alone decisive in
proving the hoped-for generalization. In our view, a key ingredient that is needed is the
phenomenon illustrated by the following example.   

\begin{example}\label{EX:MOT}
{\em Let $D = \cplx \setminus \{0,1\}$, and $f=(f_1,f_2)$ be a proper holomorphic 
  self-map of $D \times D$. Note that $D$ is hyperbolic. 
  Even though most of the hypotheses of the Remmert--Stein theorem are not satisfied, 
  the conclusion still follows.}
  
  \noindent{This is not hard to see. Fix $z_0 \in \cplx \setminus \{0,1\}$. 
  By the big Picard theorem, it follows that $0,1$ and $\infty$ are removable 
  singularities or poles of the map $h := f_1(z_0,\bcdot)$. Hence $h$ extends as a
  holomorphic map to $\widehat{\cplx}$, and is therefore a rational map. If $h$ is not
  proper as a map from $\cplx \setminus \{0,1\}$ to itself, then there is a sequence 
  $\{ x_n \} \subseteq \cplx \setminus \{0,1\}$ that converges to either $0, 1$ or 
  $\infty$, such that some subsequence of the image sequence $\{h(x_n)\}$ converges
  to a finite point in $\cplx \setminus \{0,1\}$. Hence $h$ is a rational map that 
  misses at least one of the points $0,1$ or $\infty$, and must therefore be
  constant.}

  \smallskip
  On the other hand, assume that $h$ is a proper map from $\cplx \setminus \{0,1\}$ 
  to itself; then it is a non-constant rational map. Thus, if $h=\frac{P}{Q}$, where 
  $P$ and $Q$ are two polynomials having no common factors, at least one of $P$ or
  $Q$ has to be non-constant. Also, note that $h$ takes $\{0,1,\infty\}$ to itself.

  \smallskip
  If $P$ were non-linear, it would follow that either $Q$ has the same degree as $P$,
  or $Q$ is some constant $C$. In the latter case, both $P$ and $P - Q$ are
  non-constant polynomials with disjoint zero sets. From this, it follows that $P$ is
  either $z^k$ or $(z-1)^k,\ k > 1$.  Therefore, the equation $h = 1$ has roots
  different from $0$ and $1$, which is a contradiction. If $P$ and $Q$ have the same 
  degree, then it follows that $\frac{P}{Q}$ is of the form $R^k$, where $R$ is a 
  non-constant rational function, and the value $1$ is attained by $k$ distinct 
  values, which is also a contradiction. Hence $P$ is linear, and a similar argument 
  shows that $Q$ is also linear. Hence $h$ is a fractional linear transformation that
  takes $\{0,1,\infty\}$ to itself. There are only six possibilities for the map
  $h$. From this it follows that, if for some $z_0$, $f_1(z_0,\bcdot)$ is an
  automorphism, then $f_1(z,\bcdot)$ is the same automorphism for all $z \in \cplx
  \setminus\{0,1\}$ (see Lemma \ref{LEM:CPT}). Together with the conclusion of
  the first paragraph, this proves that $f_1(z,\bcdot)$ is
  either constant for all $z$, or is independent of $z$. Applying the same argument 
  to $f_2$, we conclude that the conclusion of the Remmert--Stein theorem still
  holds.
\end{example}

The key fact used in the above example is that there are only finitely many proper
holomorphic self-maps of $\cplx \setminus \{0,1\}$. This is not true for  
the domains $\cplx$ and $\cplx \setminus \{0\}$. We are motivated by all of this
to generalize the result of Remmert and Stein to complex manifolds that are products
of certain hyperbolic Riemann surfaces. 

\begin{theorem}\label{THM:MAIN}
  Let $R_j$ and $S_j$, $j = 1, \ldots, n$, be compact Riemann surfaces, and let
  $X_j$ (resp.\,$Y_j$) be a connected, hyperbolic open subset of $R_j$
  (resp.\,$S_j$) for each $j = 1, \ldots, n$. Let $F=(F_1,\ldots,F_n) :
  X_1\times \cdots \times X_n\to Y_1\times \cdots \times Y_n$ be a finite proper
  holomorphic map. Then, denoting $z\in X_1\times \cdots \times X_n$ as
  $(z_1,\ldots,z_n)$, each $F_i$ is of the form $F_i(z_{\pi(i)})$, where $\pi$ is a
  permutation of $\{1,\ldots,n\}$.
\end{theorem}
\begin{remark}
  It is essential for $F$ to be a finite map in the above theorem. Without this
  requirement, Theorem \ref{THM:MAIN} is false. To see this, let $X$ be some
  compact hyperbolic Riemann surface. The map $F :X^2 \to X^2$ defined by
  $F(z_1,z_2) := (z_1,z_1)$ is a proper map. In fact, $F$ satisfies all the
  assumptions of Theorem \ref{THM:MAIN} except finiteness.
\end{remark}

\begin{remark}
  The conclusion of the above theorem can fail if even one of the factors is
  non-hyperbolic. Consider $X
  = \disk \times (\widehat{\cplx} \setminus \{ p \})$, where $p \in
  \widehat{\cplx}$ and $\disk$ denotes the unit disc in $\cplx$. We know that
  $\widehat{\cplx} \setminus \{ p \}$ is {\bf not} hyperbolic. We view
  $\widehat{\cplx} \setminus \{ p \}$ as $\cplx$. It is easy to check that any
  $F \in Aut(\disk \times \cplx)$ is of the form 
  \begin{equation*}
	F(z_1,z_2) = (\psi(z_1), A(z_1)z_2 + B(z_1)),
  \end{equation*}
  where $\psi \in Aut(\disk), A,B \in \hol(\disk)$ and $A$ is non-vanishing.
\end{remark}

The novelty of our proof, from the viewpoint of function theory, lies in our use of
the fact that the set of non-constant holomorphic maps between certain
Riemann surfaces is at most finite. This pheneomenon is well understood
in the realm of {\em compact} complex manifolds; see, for instance,
\cite[Chapters 6\,\&\,7]{kobayashi1998hyperbolic}. However, the factors 
$X_j$ and $Y_j$ in Theorem~\ref{THM:MAIN} are not necessarily compact. We will
see that the main idea in the Remmert--Stein theorem (i.e. Result \ref{RES:RS}) is
still useful in our more general setting. Loosely speaking, we show that, in general,
the manifold $X_1\times \cdots \times X_n$ splits into two factors, one of which 
is the product of those non-compact factors to which the Remmert--Stein method
can be applied. The finiteness result that is essential to our proof is a result by
Imayoshi \cite{imayoshi1983generalizations}. This result, plus some other technical
necessities are stated in Section~\ref{SEC:RES}. The proof of
Theorem~\ref{THM:MAIN} is presented in Section~\ref{SEC:PROOF}.

\section{A version of Montel's Theorem}
In the proof of our main result, we need a version of Montel's theorem that is
adapted to our situation. The proof of this version requires some general
results about normal families. We state these results, with references, in this
section. Throughout this section, $M$ and $N$ will denote complex manifolds, and
$\hol(M,N)$ will denote the space of holomorphic maps from $M$ into $N$. We
give $\hol(M,N)$ the compact-open topology. We begin with the definition of a normal
family.
\begin{definition}
  A subset $\mathcal{F}$ of $\hol(M,N)$ is said to be
  \emph{normal} if every sequence of $\mathcal{F}$ contains a subsequence $\{
  f_n\}$ that is either convergent in $\hol(M,N)$, or is compactly divergent. 
  By the latter we mean that given compact sets $K \subseteq M$ and $H \subseteq N$, 
  $f_n(K) \cap H = \varnothing$ for all sufficiently large $n$.
\end{definition}

\begin{result}[see \cite{kiernan1970relations}, Proposition 3]\label{RES:TAUT}
  Let $M$ be a complex manifold, and let $N$ be a complete Kobayashi hyperbolic 
  complex manifold. Then $\hol(M, N)$ is a normal family. 
\end{result}

\begin{result}[see \cite{kobayashi1967invariant}, Theorem 5.5]\label{RES:HYP}
  Let $X$ be a hyperbolic Riemann surface. Then $X$ is complete Kobayashi
  hyperbolic.
\end{result}
We now state and prove the version of Montel's Theorem that we need, which is a
corollary of the last two results. 
\begin{corollary}\label{COR:MONTEL}
  Let $X$ be a connected complex manifold and let $R$ be a hyperbolic open
  connected subset of a compact Riemann surface $S$. Then, given any sequence 
  $\{ f_{\nu}\} \subset \hol(X,R)$, there exists a subsequence  $\{ f_{\nu_k} \}$
  and a holomorphic map $f_0 : X \to \overline{R}$ (the closure taken in $S$ 
  whenever $R$ is non-compact) such that $f_{\nu_k} \to f_0$ uniformly on compact
  subsets of $X$.
\end{corollary}

\begin{proof}
  We begin by noting that if $R$ is compact, then the
  result follows immediately from Results \ref{RES:HYP} and \ref{RES:TAUT}.
  
  \smallskip
  We now consider the case when $R$ is a punctured Riemann surface. By
  Result \ref{RES:TAUT}, $\hol(X,R)$ is a normal family. 
  There is nothing to prove if there exists a subsequence
  $\{f_{\nu_k} \}$ that converges uniformly on compact subsets of $X$. Therefore,
  let us consider the case when we get only a compactly divergent subsequence 
  $\{f_{\nu_k} \}$. Let $\{K_j: j\in \mathbb{Z}_+\}$ be an exhaustion of $X$
  by {\em connected} compact subsets, and let $\{L_j: j\in \mathbb{Z}_+\}$ be
  an exhaustion of $R$ by compact subsets. Since $R$ is obtained from $S$
  by deleting finitely many points from it, compact divergence implies that
  we can extract a further subsequence from $\{f_{\nu_k} \}$ --- which we 
  shall re-index again as $\{f_{\nu_k} \}$ --- such that $f_{\nu_k}(K_1)\subset 
  D^* \; \forall k$, where $D^*$ is a deleted neighbourhood of one of the punctures, 
  say $p_0$. Now, given any $j\in \mathbb{Z}_+$, there exists a $k(j)\in 
  \mathbb{Z}_+$ such that, by the connectedness of the $K_j$'s, we have: 
  \[
   f_{\nu_k}(K_j) \subset (D^*\setminus L_j) \; \; \forall k\geq k(j).
  \]
  This just means that $f_{\nu_k} \to p_0$ uniformly on compacts as $k \to \infty$.

  \smallskip
  In the general case, as $R$ is hyperbolic, we can make sufficiently many
  punctures in $S$ to get a Riemann surface $R'$ that is hyperbolic and $R
  \subseteq R' \subset S$ . By considering each $f_{\nu}$ as a mapping in
  $\hol(X,R')$, we can find, by the preceding argument, a subsequence $\{f_{\nu_k}\}$
  and a holomorphic map $f_0:X \to \overline{R'}$ such that $f_{\nu_k} \to
  f_0$ uniformly on compact subsets of $X$. As each $f_{\nu_k} \in \hol(X,R)$,
  we must have $f_0 \in \hol(X,\overline{R})$, and we are done.
\end{proof}

\section{Some technical necessities}\label{SEC:RES}
In this section we summarize several results that we need for the proof of
Theorem \ref{THM:MAIN}. We state these results with appropriate references. 
We begin with an extension of a classical result due to de Franchis 
\cite{de1913teorema}, which states that there are at most finitely many
non-constant holomorphic mappings between two compact hyperbolic Riemann
surfaces. We shall call a Riemann surface obtained by removing a finite,
non-empty set of points from some compact Riemann surface a \emph{punctured Riemann
surface}. A Riemann surface obtained by removing $n$ points
from a compact Riemann surface of genus $g$ will be called a \emph{Riemann surface 
of finite type $(g,n)$}. Imayoshi extended de Franchis' result as follows:
\begin{result}[Imayoshi \cite{imayoshi1983generalizations}]\label{RES:FRANCHIS}
  Let $R$ be a Riemann surface of finite type and let $S$ be a Riemann surface of
  finite type $(g,n)$ with $2g-2+n > 0$. Then the set of non-constant holomorphic 
  maps from $R$ into $S$ is at most finite.
\end{result}

\noindent The above result combined with the following lemma will play a key role in 
the proof of the main theorem. To state this lemma, we need a definition. 

\begin{definition}\label{DEF:IND}
  Let $F:X \to Y$ be a map between two sets, and suppose that $X = X_1 \times \cdots
  \times X_n$. We say that $F$ is \textit{independent of $X_j$} if, for each fixed 
 $(x_1^0,\ldots,x_{j-1}^0,x_{j+1}^0,\ldots,x_n^0)$, $x_i^0 \in X_i$, the map
  \begin{equation*}
	X_j \ni x_j \longmapsto F(x_1^0,\ldots,x_{j-1}^0,x_j,x_{j+1}^0,\ldots,x_n^0),
  \end{equation*}
  is a constant map. We say that $F$ \emph{varies along} $X_j$ if $F$ is
  \textbf{not} independent of $X_j$.
\end{definition}

\begin{lemma}\label{LEM:CPT}
  Let $R$ and $S$ be as in Result \ref{RES:FRANCHIS}, and let $X$ be a connected
  complex manifold. Let $F: R \times X \to S$ be a holomorphic mapping
  with the property that for some $x_0 \in X$, the mapping $R \ni z \mapsto 
  F(z,x_0) \in S$ is a non-constant mapping. Then $F$ is independent of $X$.
\end{lemma}

\begin{proof}
  Let $d_{R}$ and $d_{S}$ be metrics that induce the topology of $R$ and $S$, 
  respectively. By Result \ref{RES:FRANCHIS}, the set of non-constant holomorphic
  mappings from $R$ to $S$ is at most finite. By our hypotheses, there is at
  least one such map. 
  Let $F_1,\ldots,F_k$ be the only distinct non-constant mappings in 
  $\hol(R,S)$. Let $x_0 \in X$ be such that the map $F(\bcdot,x_0)$ is
  non-constant. By continuity of $F$, there is an $X$-open neighbourhood 
  $U_0 \ni x_0$ such that $F(\bcdot,x)$ is non-constant for $x \in U_0$. 
  Choose $\eps > 0$ and $r_{ij} \in R,\ 1 \leq i,j  \leq k, \ i \neq j$, such that 
  $d_{S}(F_i(r_{ij}),F_j(r_{ij})) > \eps$. By the continuity of $F$, we can find 
  a neighbourhood $U \subset U_0$ of $x_0$ such that, for each of the $r_{ij}$'s, we
  have $d_{S}(F(r_{ij},x),F(r_{ij},y)) < \eps \ \forall x,y \in U$. This is possible
  only if $F(\bcdot,x) \equiv F(\bcdot,y),\ \forall x,y \in U$. It follows that
  \begin{itemize}
	\item $\exists j_0 \leq k$ such that $F(\bcdot,x)=F_{j_0} \ \forall x \in U$;
	\item For each fixed $r \in R$, the map $F(r,\bcdot)$ is constant on $U$.
  \end{itemize}
  As $X$ is connected, the Identity Theorem implies that $F(r,\bcdot) \equiv
  F_{j_0}(r)$. This proves that $F$ is independent of $X$.
\end{proof}

The next result is the well known Remmert's Proper Mapping Theorem. For the proof,
refer to \cite[p.~31]{chirka1989complex}. 

\begin{result}[Proper Mapping Theorem]\label{RES:PROPER}
Let $X$ and $Y$ be complex manifolds, and let $A$ be an analytic subset of $X$.
Let $f: A \to Y$ be a proper finite holomorphic map. Then, $f(A)$ is an analytic
subset of $Y$, and at every $w \in f(A)$
\begin{equation*}
  dim_w f(A) = \max \{dim_z A : f(z) = w \}.
\end{equation*}
In particular, $dim A = dim f (A)$. Furthermore, if $A = X$ and $dim(X) = dim(Y)$ 
then $F$ is surjective.
\end{result}

The following result due to Kobayashi \cite[p.~284]{kobayashi1998hyperbolic}
can be thought of as an higher dimensional analogue of the big Picard theorem.
For this, we first need to make a definition.

\begin{definition}
  Let $Z$ be a complex manifold and let $Y$ be a relatively compact complex 
  submanifold of $Z$. We call a point $p \in \overline{Y}$ a 
  \textit{hyperbolic point} if every neighbourhood $U$ of $p$ contains a smaller
  neighbourhood $V$ of $p,\ \overline{V} \subset U$, such that
  \[
    K_{Y}(\overline{V} \cap Y, Y\setminus U) := \inf\{K_Y(x,y) : x \in
	  \overline{V} \cap Y, y \in Y \setminus U\} > 0,
  \]
  where $K_{Y}$ denotes the Kobayashi pseudo-distance on $Y$. We say that $Y$ is
  \emph{hyperbolically imbedded} in $Z$ if every point of $\overline{Y}$ is a
  hyperbolic point.  
\end{definition}

\begin{result}[Kobayashi]\label{RES:BPIC}
  Let  $Y$ and $Z$ be complex manifolds, and let $Y$ be hyperbolically imbedded 
  in $Z$. Then every map $h \in \hol\left(\disk ^*,Y\right)$ extends to a map 
  $\tilde{h} \in \hol(\disk,Z)$.
\end{result}

\begin{lemma}
  If $Y$ is a hyperbolic open connected subset of a compact Riemann surface $Z$, 
  then $Y$ is hyperbolically imbedded in $Z$.
\end{lemma}

\begin{proof}
 The lemma is obvious if $Y$ has only isolated boundary points. If not, then, as
 $Y$ is hyperbolic, we can make sufficiently many punctures in $Z$ to get a
 hyperbolic Riemann surface $\widetilde{Y}$ such that $Y \subset \widetilde{Y}
 \subset Z$. It follows that $\widetilde{Y}$ is hyperbolically imbedded in $Z$.
 Now let $y \in \overline{Y}$. Then, $y \in \overline{\widetilde{Y}}$. Let
 $U$ be a neighbourhood of $y$, and let $V$ be a smaller neighbourhood of $y$
 such that
 \[
   K_{\widetilde{Y}}(\overline{V} \cap \widetilde{Y}, \widetilde{Y}\setminus U)
   > 0.
 \]
 As $K_Y \geq K_{\widetilde{Y}}$ on $Y\times Y$ and $Y \subseteq \widetilde{Y}$, it follows 
 that $y$ is also a hyperbolic point of $Y$. Consequently, $Y$ is hyperbolically 
 imbedded in $Z$.
\end{proof}

\noindent We require one more result, a generalization of Rado's theorem that is
proved in \cite{narasimhan1971scv}.

\begin{lemma}\label{RES:RADO}
  Let $(\phi_{\mu\nu}),\ 1 \leq \mu \leq k,\ 1 \leq \nu \leq l$, be a matrix of
  holomorphic functions on $D \subseteq U$, where $D$ and $U$ are connected open
  subsets of $\cplx$, and $U\setminus D$ is a non-empty indiscrete set. 
  Suppose that
  \[
	\prod_{\nu=1}^l\sum_{\mu=1}^k\left|\phi_{\mu\nu}(z)\right|^2 \to 0
    \ \text{as} \  D \ni z \to \zeta
  \]
  for any $\zeta \in \partial D \cap U$. Then, for some $\nu_0, \ 1\leq \nu_0
  \leq l$, we have 
  \[
  \phi_{\mu\nu_0} \equiv 0, \ \mu=1,\ldots,k.
  \]
\end{lemma}

\begin{proof}
  Suppose each column of $(\phi_{\mu\nu})$ has a member that is not
  identically $0$ on $D$. Let $f$ be the product of these members. We extend $f$
  to be a function on $U$ by defining $f \equiv 0$ on $U \setminus D$. By 
  hypothesis, $f$ is continuous on $U$ and holomorphic on $D$. Therefore by the 
  classical Rado's theorem $f \equiv 0$, a contradiction. 
\end{proof}

\section{Proof of the main theorem} \label{SEC:PROOF}

We begin this section by considering a special case of Theorem \ref{THM:MAIN}
whose proof contains some technicalities. Since these technicalities would lengthen
the proof of Theorem \ref{THM:MAIN} if we were to embark on it 
directly, we shall isolate the technical portion of our proof in the following
proposition. Its proof consists of rephrasing the Remmert--Stein argument
relative to a coordinate patch; see \cite[pp.~71--78]{narasimhan1971scv}.
{\em We shall therefore be brief} and explain in detail only those points that differ
from the proof in \cite{narasimhan1971scv}.
\begin{proposition}\label{PROP:BORD}
  Let $X=X_1\times \cdots \times X_n$ and $Y=Y_1 \times \cdots \times
  Y_n,\ n \geq 2$, be complex manifolds. Assume that each $X_j$  and each $Y_j$
  satisfy the hypothesis of Theorem \ref{THM:MAIN} and that $Y$ is non-compact.
  Further assume that, for each $j$, $R_j \setminus X_j$ is a non-empty indiscrete set. 
  Let $F: X \to Y$ be a finite proper holomorphic map. Then, denoting $z\in X_1
  \times \cdots \times X_n$ as $(z_1,\ldots,z_n)$, each $F_i$ is of the form
  $F_i(z_{\pi(i)})$, where $\pi$ is a permutation of $\{1,\ldots,n\}$. 
  
  \smallskip
  \noindent In particular, if there is a mapping with the
  above properties from $X$ to $Y$, then $Y$ cannot have any compact factors.
\end{proposition}
\begin{proof}
  For $1 \leq j \leq n$, let $R_j$ and $S_j$ be as in Theorem \ref{THM:MAIN}. Let 
  $p=(p_1,\ldots,p_n)$ be a point in $R_1 \times \cdots \times R_n$ such
  that, for each $1 \leq i \leq n,\ p_i$ is a limit point of the set $R_i
  \setminus X_i$ and belongs to $\bdy X_i$.

  Let $(U_k,\psi_k)$ be connected holomorphic co-ordinate charts of $R_k$, chosen in
  such a way that $p_i \in U_i$ and the image of $(\prod_{k=1}^n U_k) \cap X$ under 
  each $F_j$ lies in some holomorphic co-ordinate chart $(V_j,\rho_j)$ of $S_j$. Let
  $W_i$ be a connected component of $U_i \cap X_i$ such that $p_i \in \bdy W_i$.
  For $(z_1,\ldots,z_n) \in \prod_{k=1}^n \psi_k(W_k)$, let
  \begin{equation*}
	g_j(z_1,\ldots,z_n) : =  \rho_j\circ
	F_j(\psi_1^{-1}(z_1),\ldots,\psi_n^{-1}(z_n)).	
  \end{equation*}
   In view of Corollary \ref{COR:MONTEL}, we can rephrase the
   arguments in \cite[p.~75]{narasimhan1971scv} to conclude:
  \[
	\prod_{j=1}^n \sum_{k=1,k \neq i}^n \left|\frac{\partial g_j}{\partial
	z_k}(z_1,\ldots,z_{i-1},w,\ldots,z_n) \right|^2 \to 0 \; \ \text{as} \; \ 
	w \to \zt\in \psi_i(\partial W_i),
  \]
  where $\zt$ is any arbitrary point in $\psi_i(\partial W_i)$.
  Let us take $D = \psi_i(W_i)$ and $U = \psi_i(U_i)$ in Lemma \ref{RES:RADO}.
  Note that $\psi_i(p_i) \in U \setminus D$, whence $U \setminus D$ is
  indiscrete. Thus, we have that 
  for each $(z_1,\ldots,\ z_{i-1},\ z_{i+1},\ldots,\ z_n) \in 
  \prod_{k = 1, k \neq i}^n \psi_k(W_k)$, there is a $j = j(z)$ such that 
  \begin{equation*}
	h_j(z_1,\ldots,z_{i-1},w,z_{i+1},\ldots,z_n) := 
	\sum_{k=1,k \neq i}^n \left|\frac{\partial g_j}{\partial
	z_k}(z_1,\ldots,z_{i-1},w,\ldots,z_n) \right|^2
  \end{equation*}
  is zero $\forall w \in D$. At this point, we can again argue exactly as in
  \cite[p.~75]{narasimhan1971scv} to conclude that there exists an integer
  $\sigma(i), 1 \leq \sigma(i) \leq n$, such that
  \[
  \frac{\partial g_{\sigma(i)}}{\partial z_k} \equiv 0 \; \ \text{on} \; \ \psi_1(W_1) 
  \times \cdots \times \psi_n(W_n), \; \ k = 1,\ldots,n, \; \ k \neq i.
  \]
  Therefore on $W_1 \times \cdots \times W_n ,\ F_{\sigma(i)}$ is independent of
  $z_1,\ldots,z_{i-1},z_{i+1},\ldots,z_n$. By applying the Identity Theorem, 
  we conclude that $F_{\sigma(i)}$ is independent of the same variables on $X$. By
  Remmert's Proper Mapping Theorem, $F$ is surjective. This implies that 
  $F_{\sigma(i)}$ varies along $X_{i}$. Since the choice of $1 \leq i \leq n$ in the
  preceding argument was arbitrary, for each $i$ there exists precisely one
  $\sigma(i)$ such that $F_{\sigma(i)}(z) =
  F_{\sigma(i)}(z_i) \ \forall z \in X$. The permutation $\pi = \sigma^{-1}$, and
  we are done with the proof of the first part.
  
  \smallskip
  To establish the final part of this result, assume that $Y_{s+1},\ldots,Y_n$
  are all compact, for some $s < n$. Fix an $i$ as in the previous paragraph.
   The heart of the argument above, see
  \cite[p.~75]{narasimhan1971scv}, consists of using Montel's theorem 
  (Corollary~\ref{COR:MONTEL} in our present set-up) to construct a map
  $(\phi_1,\ldots,\phi_n) : Z \to \partial Y$, where $Z:= \prod_{k=1,k \neq i}^n X_k$.
   Set $E_j := \{z\in Z: \phi_j(z)\in \partial Y_j\}$. Clearly :
  \begin{equation} \label{EQ:1}
	\{ l : 1 \leq l \leq n, \ \text{int}(E_l) \neq \varnothing \} \subseteq \{
	1, \ldots,s\}.
  \end{equation}
  In view of (\ref{EQ:1}), the argument in \cite[p.~75]{narasimhan1971scv} reveals
  that, for each $i, \sigma(i) \in \{ 1,\ldots,s \}$. Since $s < n$, by assumption, there
  would exist $i \neq i'$ such that $\sigma (i) = \sigma (i')$. But this would
  contradict the surjectivity of $F$, and we are done.

\end{proof}

\begin{custom}\label{CUS:PROOF}{\bf The proof of Theorem \ref{THM:MAIN}.} 
  For $1 \leq j \leq n$, let  $R_j$ and $S_j$ be the compact 
  Riemann surfaces in the statement of the theorem.
  We start off with a simple consequence of the finiteness of $F$.\\
  \textbf{Claim A}: \emph{For any holomorphic finite map $F:X \to Y$, given any $X_i, 
  1 \leq i \leq n$, there is some $F_j$ that varies along $X_i$.} 
  To see this, assume that there is a factor $X_i$ such that all the $F_j$'s are
  independent of $X_i$. Then for any point $x = (x_1,\ldots,x_n) \in X$, by 
  Definition \ref{DEF:IND}, the inverse image of $F(x)$ contains the set
  $\{x_1\} \times \cdots \{x_{i-1}\} \times X_i \times \{x_{i+1} \} \times \cdots
  \{ {x_n}\}$. But this contradicts the finiteness of $F$. 

  \smallskip
  Let $X_C$ and $Y_C$ denote the product of those factors of $X$ and $Y$, 
  respectively, that are either compact, or compact with finitely many punctures, 
  and let $X_B$ and $Y_B$ denote the product of the remaining factors. Since 
  Proposition \ref{PROP:BORD} already establishes our theorem if $X_C = \varnothing$,
  we may assume, without loss of generality, that $X_C := X_1 \times \ldots \times
  X_p, \ 1 \leq p \leq n$. Note that if $X_C = \varnothing$ and
  there exists a proper holomorphic map $F: X \to Y$, then $Y$ cannot be compact.

  \smallskip
  \noindent \textbf{Claim B}: \emph{The maps
  $F_i$ are independent of $X_1,\ldots, X_p$, whenever $S_i \setminus Y_i$ is
  a non-empty indiscrete set.} 
  To see this, fix $x' \in X_2 \times \cdots \times X_n$. The map 
  $F_i(\bcdot,x')$ is a holomorphic map from $X_1$ into a hyperbolically imbedded
  Riemann surface. Now, $R_1$ is the compact Riemann surface from which $X_1$ is 
  obtained by deleting at most finitely many points. From Result \ref{RES:BPIC}, it 
  follows that $F_i(\bcdot,x')$ extends holomorphically to a map $\widetilde{f}_i$
  from $R_1$ into $S_i$. If $\widetilde{f}_i$ is non-constant, then by the 
  compactness of $R_1$, it follows that the image set of $R_1$ under 
  $\widetilde{f}_i$ is both compact and open. But, this means that $S_i = 
  \widetilde{f}_i(R_1)$, which is not possible as $S_{i} \setminus Y_i$
  is a non-empty indiscrete set, and $\widetilde{f}_i(R_1)$ is obtained by adjoining at
  most finitely many points to $Y_i$. This proves that $F_i$ is independent of $X_1$. 
  Repeating the same argument for the factors $X_2,\ldots,X_p$, the claim is proved.

  \smallskip
  Now note that, in view of Claim B, if $1 \leq i \leq p$ and $F_j$ is a
  map that varies along $X_i$, then $Y_j$ is either compact, or compact with finitely
  many punctures. Then,
  by Lemma \ref{LEM:CPT}, $F_j$ is independent of \emph{all the factors of $X$ other 
  than $X_i$.} Without loss of generality, we may assume that $Y_C = Y_1
  \times \ldots Y_k, \ 1 \leq k \leq n$. Combining our last deduction with Claim
  A, we infer that:
  \begin{enumerate}
	\item $p \leq k \leq n$;
	\item Without loss of generality, there is an enumeration of the factors of
	  $Y_C$ such that for each $1 \leq i \leq p$, there is a unique $\sigma(i),
	  \ 1 \leq \sigma(i) \leq p$, such that $F_{\sigma(i)}(z) =
	  F_{\sigma(i)}(z_i) \ \forall z \in X$.
  \end{enumerate}

  \smallskip
  Suppose $k > p$. Then, in view of the (harmless) assumption in (2), we need
  to analyse the behaviour of $F_i, \ p+1 \leq i \leq k$. Note that we already
  know from Claim B that $F_{k+1},\ldots,F_n$ is independent of $X_C$. 
  Assume that $F_{p+1}$ varies along some $X_i, 1 \leq i \leq p$; then from Lemma 
  \ref{LEM:CPT}, $F_{p+1}$ is independent of all other factors of $X$. From
  Remmert's Proper Mapping Theorem (Result \ref{RES:PROPER}), $F$ is a surjective
  map from $X$ onto $Y$. Hence, combining the last two assertions with (2),
  $(F_1,\ldots,F_{p+1})$ determines a surjective holomorphic
  map $(F_1,\ldots,F_{p+1}): X_C \to Y_1 \times \cdots \times Y_{p+1}$ from a space
  of dimension $p$ to a space of dimension $p+1$, which contradicts Sard's theorem. 
  Hence, $F_{p+1}$ is independent of $X_1,\ldots,X_p$.
  Repeating the same argument for each map
  $F_j, \ p+1 \leq j \leq k$, we conclude that each $F_j, \ p+1 \leq j \leq n$, is
  independent of $X_C$.

  \smallskip
  Whether or not $k > p$, the previous paragraph implies that $F_i, \ p+1 \leq i 
  \leq n$, are independent of $X_C$, whence they determine a
  surjective map $F_B = (F_{p+1},\ldots,F_n) : X_B \to Y_{p+1} \times \ldots
  \times Y_n$. This map is clearly finite. We will now show that it is proper.
  Consider a compact set $K \subseteq  Y_{p+1} \times \ldots \times Y_n$. We
  must show that $F_B^{-1}(K)$ is a compact subset of $X_B$. Let $H \subseteq Y_1
  \times \cdots \times Y_p$ be some compact set. Then, by the properness of $F$,
  it follows that $F^{-1}(H \times K)$ is compact. But, given the independence
  of the various $F_i$'s from certain factors of $X$, 
  \begin{equation*}
	F^{-1}(H \times K) = (F_1,\ldots,F_p)^{-1}(H) \times F_B^{-1}(K).	
  \end{equation*}
  Thus, $F_B^{-1}(K)$ is compact, as required.
  
  \smallskip
  As $X_B$ is non-compact, and $F_B$ is a proper map, it follows that $Y_{p+1}
  \times \cdots \times Y_n$ is also non-compact. We now apply Proposition
  \ref{PROP:BORD} to the map $F_B$ to get a permutation $\pi$ of $\{p+1,\ldots,n\}$ 
  such that, for each $p < i \leq n$, we have $F_i(z) = F_i(z_{\pi(i)})$.
  Juxtaposing $\pi$ with the permutation $\sigma^{-1}$ of $ \{1,\ldots,p\}$, we are
  done.
  \qed
\end{custom}

\noindent {\bf Acknowledgement.} I would like to thank my advisor Gautam Bharali
for his support during the course of this work, and for suggesting several
useful ideas. I would also like to thank my colleagues and friends G.P.
Balakumar, Dheeraj Kulkarni and Divakaran Divakaran for the many interesting
discussions.

\bibliographystyle{amsplain}
\bibliography{hyperbolic}

\end{document}